\documentstyle[11pt]{article}

\bibstyle{plain}

\newcommand{\al}{\alpha}
\newcommand{\R}{{\mathbf R}}
\newcommand{\ds}{\displaystyle}
\newcommand{\e}{\varepsilon}
\newcommand{\Om}{\Omega}
\newcommand{\lra}{\longrightarrow}
\newcommand{\ra}{\rightarrow}
\newcommand{\p}{\partial}
\newcommand{\la}{\lambda}
\newcommand{\g}{\gamma}

\newcommand{\Z}{{\mathbf Z}}
\newcommand{\Q}{{\mathbf Q}}

\newcommand{\te}{\theta}

\newcommand{\ph}{\varphi}

\newcommand{\Po}{{\mathcal P}}
\newcommand{\Xo}{{\mathcal X}}
\newcommand{\Bo}{{\mathcal B}}
\newcommand{\Mo}{{\mathcal M}}
\newcommand{\DR}{{\mathcal D}^{1,2}(\R^N)}

\newcommand{\wh}{\widehat}

\newcommand{\references}[1]{\theinstitutions 
\begin{thebibliography}{#1}}

\newcommand{\Endrefs}{\end{thebibliography}}

\title{\bf On the symmetry of minimizers}

\date{  }

\author{ Mihai MARI\c S \\
          { \small     } \\
         \it  Universit\'e de Franche-Comt\'e \\
         \it   D\'epartement de Math\'ematiques UMR 6623\\
         \it   16, Route de Gray \\
         \it 25030 Besan\c con, France \\
         \it  e-mail: mihai.maris@univ-fcomte.fr}

\date{}

\begin{document}

\maketitle

\vspace{-20pt}
\begin{center}
\it Dedicated to  Dorel Mihe\c t, \\
for his teaching, his friendship, and the inspiration he gave to me.
\end{center}
\rm

\vspace{-5pt}

\noindent
\begin{abstract}
For a large class of variational problems we prove that minimizers are symmetric whenever they are $C^1$.

\medskip

\noindent
{\bf AMS subject classifications. } 35A15, 35B05, 35B65, 35H30, 35J20, 35J45, 35J50, 35J60.

\end{abstract}

\section{Introduction and main results}

In this paper we study the symmetry of minimizers for general variational problems of the form
$$
\begin{array}{c}
\mbox{ minimize } E(u) := \ds \int_{\Om} F(|x|, u(x), |\nabla u(x)|) \, dx 
\qquad \mbox{ under $k$ constraints }
\\
\\
Q_j(u) = \ds \int_{\Om} G_j(|x|, u(x), |\nabla u(x)|) \, dx = \la _j, 
\qquad j=1, \dots, k.
\end{array}
\leqno{(\mathcal{P})}
$$

The solutions of many  partial differential equations are obtained as 
minimizers for problems like ($\Po$).
Knowing in advance that such solutions are symmetric is very important 
for their theoretical study as well as for their numerical approximation.
If the minimizers of ($\Po$) are standing or solitary waves for 
an evolution equation, symmetry could be very useful to investigate 
the stability properties of such solutions.
Note also that in many problems, symmetry is the first step in proving the uniqueness of 
special solutions.

Given the motivation above, many important particular cases of ($\Po$) 
have already been considered in the literature.
In \cite{lop1, lop2}, O. Lopes has developed his reflection method - a very efficient 
tool to prove symmetries for minimizers of functionals 
$E_1(u) = {\ds \int_{\Omega}} \frac 12 |\nabla u |^2 + F_1(|x|, u) \, dx $
under the constraint  
$Q (u) = \ds \int_{\Omega} G(|x|, u) \, dx = constant, $
where $ \Om $ is a domain invariant by rotations. 
This method is based on a device of "reflecting" a minimizer with respect to 
hyperplanes that "split the constraint in two" and on the use of a unique 
continuation principle for the Euler-Lagrange equations satisfied by  minimizers.
Note that 
the method can be used for vector-valued minimizers 
whose components eventually change sign and no additional assumptions are made on the functions $F_1$ and $G$ (except the usual growth and
smoothness assumptions that ensure the existence and the regularity of minimizers).
Up to now this method has been used for problems involving only one constraint.
Its main restriction is that it can be used only when the minimizers satisfy an 
Euler-Lagrange system for which a unique continuation theorem is available.
However, we have to mention that the reflection 
method has been successfully used in \cite{LM}
for minimizers of some  nonlocal functionals of the form 
$E_2(u) = {\ds \int_{\R^N}} m(\xi ) |\wh{u}(\xi ) |^2\, d \xi + 
{\ds \int_{\R^N}} F_2 (u) \, dx$. 
The class of functionals considered in 
\cite{LM} include 
the generalized Choquard functional, the Hamiltonian for the generalized 
Davey-Stewartson equation as well as functionals involving 
fractional powers of the Laplacian.
Instead of unique continuation results, some new and quite unexpected 
integral identities for nonlocal operators were used 
to get symmetry results.

In a recent paper \cite{brock-JMAA}, 
F. Brock studies the symmetry of minimizers of the functional
${\ds \int_{\R^N}} \ds \sum_{i =1}^n |\nabla u_i|^p + F(|x|, u_1, \dots, u_n)\, dx $ under several constraints 
${\ds \int_{\R^N}} G_{i,j} (u_i)\, dx = c_{i,j}$.
He uses two-points rearrangements 
and a variant of the strong maximum principle due to Pucci, Serrin and Zou (\cite{PSZ})
to prove symmetries.
Assuming that $F$ is nonincreasing in the first variable and that 
$\frac{\p F}{\p u_i}$ is nondecreasing in the variables $u_k$ for $ k \neq i$
(a cooperative condition), 
he shows that the superlevel sets $\{u_i >t \}$ for $t >0$,
respectively the sublevel sets $\{u_i <t \}$ for $t <0$,  are balls.
Under more restrictive conditions 
($F$ strictly decreaing in the first variable or an assumption that  depends
 on  Lagrange multipliers associated to minimizers - 
assumption that  could be quite difficult to check in applications, as already mentioned in \cite{brock-JMAA}), he proves that any component of the minimizer is radially symmetric about $0$, 
has constant sign and is monotone in $|x|$. 
Note that whenever the arguments in \cite{brock-JMAA} lead to symmetry, 
they also imply monotonicity.
On the other hand, in \cite{brock-JMAA} there is an example of 
sign-changing minimizer for a particular functional of the type considered. 
It is remarkable that the results of F. Brock are valid for an arbitrary number of constraints.
However, these constraints must have a special form (because they have to be preserved by rearrangements of functions). 
For instance, one cannot allow constraints of the form 
${\ds \int_{\R^N}} G(u_i, u_j)\, dx = constant$.

We have to mention that in a series of recent papers (see \cite{BWW}, \cite{PW}, \cite{SW}
and references therein), different new techniques were developed to study 
the symmetry of solutions for some classes of elliptic problems. 
These techniques are essentially based on foliated Schwarz rearrangements 
and on polarization of functions and can be used for sign-changing solutions. 
They also give some monotonicity properties. 

\medskip

The aim of the present paper is to prove symmetry of minimizers for problem 
($\Po$) under general assumptions.
We use the device of reflecting minimizers with respect to hyperplanes introduced by O. Lopes, but we do not need unique continuation theorems. 
Instead, we use in an essential way the regularity of minimizers.
(To our knowledge,   symmetry results  for 
 minimizers that may be nonsmooth  were obtained 
only in the case of convex functionals.)  
We are able to deal with several constraints, but each additional constraint produces 
the loss of one direction of symmetry;  we will see later (Examples 6 and 7)
that under the general conditions  considered here, 
this is a very natural phenomenon.

\medskip

In the sequel $\Om $ denotes an open set in $\R^N$ invariant under rotations 
(and centered at the origin).
It is not assumed that $\Om $ is connected or bounded.
We denote $A_{\Om} = \{ |x| \; | \; x \in \Om \}$.
We  consider vector-valued minimizers $ u : \Om \lra \R^n$ 
of ($\Po$) that belong to some function space $\Xo$. 
Throughout $F, G_1, \dots G_k$ are real-valued functions defined on 
$ A_{\Om} \times \R^m \times [0, \infty )$ in such a way that 
for any $ v \in \Xo$, 
the functions $ x \longmapsto F(|x|, v(x), |\nabla v (x)|) $ and 
$ x \longmapsto G_j(|x|, v(x), |\nabla v (x)|) $, $ 1 \leq j \leq k$, 
belong to $L^1(\Om )$. 

\medskip

Let $V$ be an affine subspace of $\R^N$. 
For $ x \in \R^N$ we denote by $p_V (x) $ the projection of $x$ onto $V$ and 
by $s_V(x)$ the symmetric point of $x$ with respect to $V$,  
that is $ s_V(x) = 2p_V(x) -x$.
We say that a function $f$ defined on $\R^N$ is 
{\it symmetric } with respect to $V$ if $f(x) = f(s_V(x)) $ for any $x$. 
We say that $f$ is {\it radially symmetric } with respect to $V$ if 
there exists a function $\tilde{f}$ defined on $ V \times [0, \infty )$ 
such that 
$f(x) = \tilde{f} (p_V(x), |x - p_V(x)|)$.

\medskip

Let $\Pi $ be a hyperplane in $\R^N$ and let $\Pi^+$ and $ \Pi ^-$ be the two half-spaces determined by $\Pi$.
Given  a function $f$ defined on $\R^N$, we denote 
$$
\begin{array}{l}
f_{\Pi^+}(x) = \left\{ 
\begin{array}{lcl}
f(x) & \mbox{ if } & x \in \Pi^+ \cup \Pi , \\
f(s_{\Pi} (x)) & \mbox{ if } & x \in \Pi^- , 
\end{array}
\right. 
\qquad \mbox{ respectively }  \\ \\
f_{\Pi^-}(x) = \left\{ 
\begin{array}{lcl}
f(x) & \mbox{ if } & x \in \Pi^- \cup \Pi , \\
f(s_{\Pi} (x)) & \mbox{ if } & x \in \Pi^+ . 
\end{array}
\right. 
\end{array}
\leqno{(1)}
$$
If $f$ is defined  on 
a rotation invariant subset
$ \Om  $ centered at the origin,  $\Om \neq \R^N$,  
the above definition makes sense only if 
$\Pi$ contains the origin.
We say that $\Pi $ splits the constraints in two for a function $ v \in \Xo$ 
if 
$$
\ds \int_{\Om \cap \Pi ^+} G_j(|x|, v(x) , |\nabla v(x)|) \, dx = 
\ds \int_{\Om \cap \Pi ^-} G_j(|x|, v(x) , |\nabla v(x)|) \, dx
\quad 
\mbox{ for } j=1, \dots, k.
\leqno{(2)}
$$

We make the following assumptions. 

\bigskip

\noindent
{\bf A1. } For any $ v \in \Xo $ and any hyperplane $\Pi$ containing the origin, we have $ v_{\Pi ^+}, v_{\Pi ^-} \in \Xo$. 

\medskip

\noindent
{\bf A2. } Problem ($\Po$) admits minimizers in $\Xo$ and any  minimizer is a $C^1$ function on $\Om$. 

\bigskip

We can now state our symmetry results.

\bigskip

\noindent
{\bf Theorem 1. } \it Assume that $ 0 \leq k \leq N-2$ and {\bf A1, A2 } are satisfied. 
Let $ u \in \Xo $ be a minimizer for problem ($\Po$). 
There exists a $k-$dimensional  vector subspace $V$ in $\R^N$ such that 
$u$ is radially symmetric with respect to $V$. 

\rm 

\bigskip

If $ \Om = \R^N$ and the considered functionals are invariant by translations, Theorem 1 can be improved. 
More precisely, consider the following particular case of ($\Po$):
$$
\begin{array}{c}
\mbox{ minimize } E(u) := \ds \int_{\R^N} F( u(x), |\nabla u(x)|) \, dx 
\qquad \mbox{ subject to  $k$ constraints }
\\
\\
Q_j(u) = \ds \int_{\R^N} G_j( u(x), |\nabla u(x)|) \, dx = \la _j, 
\qquad j=1, \dots, k.
\end{array}
\leqno{(\Po ')}
$$
In this case  assumption {\bf A1} is replaced  by 

\bigskip

\noindent
{\bf A1.' } For any $ v \in \Xo $ and any affine hyperplane $\Pi$ 
in $\R^N$ 
 we have $ v_{\Pi ^+}, v_{\Pi ^-} \in \Xo$. 

\bigskip

The following result holds.

\bigskip

\noindent
{\bf Theorem 2. } \it Assume that $ 1 \leq k \leq N-1$, {\bf A1'}
and {\bf  A2 } are satisfied and there exists $ j \in \{ 1, \dots, k \}$ such that 
$ \la _j \neq 0$. 
Let $ u \in \Xo $ be a minimizer for problem ($\Po '$). 
There exists a $(k-1)-$dimensional  affine subspace $V$ in $\R^N$ such that 
$u$ is radially symmetric with respect to $V$. 

\rm 

\medskip

If ($\Po '$) involves only one constraint, Theorem 2 implies that
any minimizer is radial with respect to some point.

\medskip

In applications, assumptions {\bf A1} or {\bf A1'} are usually  easy to check. 
On the contrary, assumption {\bf A2} requires much more attention.
In most applications, under suitable growth and smoothness assumptions on the 
functions $F, G_1, \dots , G_k$, 
the functionals $E, Q_1, \dots , Q_k$ are differentiable on 
$\Xo $ and the minimizers satisfy Euler-Lagrange equations
(however, this is not always the case: see \cite{ball-mizel} 
for examples of minimizers that do not satisfy Euler-Lagrange equations). 
Very often the Euler-Lagrange equations are, in fact, quasilinear elliptic systems. 
Many efforts have been made during the last 50 years, since the pioneer work of 
de Giorgi, Nash and Moser, to study the regularity of solutions of such systems and there is a huge literature devoted to the subject. 
Important progress has been made and various sufficient conditions that guarantee the regularity of solutions have been given. 
It would exceed the scope of the present paper to resume these works, 
or even to give here a significant list of conditions that ensure the regularity of minimizers.
For these issues (and also for historical notes) we refer the reader to 
the standard books \cite{chen-wu, giaquinta1, giaquinta2, GT, LU, maly-ziemer} and references therein.

In the next section we give the proofs of Theorems 1 and 2. 
We end this paper by some remarks and examples which show that, 
under the general conditions considered here, our results are optimal
even for scalar-valued minimizers.

\section{Proofs }

$\; $

\vspace{-11pt}

{\it Proof of Theorem 1. }
 Consider first the case $ 1 \leq k \leq N-2$. 
For $ {v} \in \R^N$, $v \neq 0$, denote 
$ \Pi_{v } = \{ x \in \R^N \; | \; x. v =0 \}$, 
$ \Pi_{v } ^+  = \{ x \in \R^N \; | \; x. v > 0\}$ and
$ \Pi_v ^-= \{ x \in \R^N \; | \; x.v < 0\}$.
For $ j = 1, \dots, k$, we define $  \ph _j : S^{N-1} \lra  \R$ by 
$$
\ph _j(v) = 
\ds \int_{\Pi_{v }^+ \cap \Om } G_j(|x|, u(x), |\nabla u(x)|) \, dx 
- \ds \int_{\Pi_ {v }^-\cap \Om } G_j(|x|, u(x), |\nabla u(x)|) \, dx .
$$
It is obvious that $\ph_j(-v ) = \ph_j (v)$ and it follows immediately from 
Lebesgue's dominated convergence theorem that each $\ph _j$ is continuous on 
$S^{N-1}$. We will  use  the following well-known result
(see, e.g., \cite{spanier}, Theorem 9 p. 266):

\bigskip

\noindent
{\bf Borsuk-Ulam Theorem. } \it Given a continuous map 
$ f : S^{n_1} \lra \R^{n_2}$ 
with $ n_1 \geq n_2 \geq 1$, there exists $ x \in S^{n_1}$ such that
 $ f(x) = f(-x)$. 

\rm   

\bigskip

Equivalently, any continuous odd map $ f : S^{n_1} \lra \R^{n_2}$, 
$ n _1 \geq n_2 \geq 1$, must vanish. 

We use the Borsuk-Ulam theorem for the odd continuous map 
$ \Phi = ( \ph_1, \dots , \ph _k) : S^{N-1} \lra \R^k$ and we infer that 
there exists $ e _1 \in S^{N-1}$ such that $ \Phi(e_1) = 0$, 
that is 
$ \Pi_{e_1} $ splits the constraints in two for the minimizer $u$. 

Our aim is to show that $ u$ is symmetric with respect to $ \Pi_{e_1} $.
We denote $u_1 = u_{\Pi_{e_1}^-} $ and  $u_2 = u_{\Pi_{e_1}^+} $
the two reflected functions obtained from $u$ as in (1). 
By {\bf A1 } we have $ u_1, \; u_2 \in \Xo$. 
Since $\Pi_{e_1} $ splits the constraints in two, a simple change of variables shows that 
$
\ds \int_{\Om } G_j(|x|, u_1(x), |\nabla u_1 (x)|) \, dx 
= 2 \ds \int_{\Pi_v^- \cap \Om } G_j(|x|, u_1 (x), |\nabla u_1 (x)|) \, dx 
= \la _j $ for any $ j \in \{ 1, \dots, k\}$, that is $ u_1 $ satisfies the constraints. 
In the same way $ u_2$ satisfies the constraints.
Since $u$ is a minimizer for ($\Po$), we must have 
$ E(u_1) \geq E(u) $ and $ E(u_2) \geq E(u) $.
On the other hand, we get 
$$
\begin{array}{l}
E(u_1) + E(u_2)  = 
2\ds \int_{\Pi_v^- \cap \Om } F(|x|, u_1 (x), |\nabla u_1 (x)|) \, dx 
+ 
2\ds \int_{\Pi_v^+\cap \Om } F(|x|, u_1 (x), |\nabla u_1 (x)|) \, dx 
\\
 =  2E(u).
\end{array}
$$
Thus necessarily $ E(u_1) = E(u_2) = E(u)$ and $ u_1$, $u_2$ are also minimizers for problem ($\Po$).
Moreover, they are symmetric with respect to $ \Pi_{e_1}$.

Now let us consider the minimizer $ u_1$. 
We define $ \psi _j : S^{N-1} \lra \R$ by 
$$
\psi _j(v) = 
\ds \int_{\Pi_v^+ \cap \Om } G_j(|x|, u_1(x), |\nabla u_1(x)|) \, dx 
- \ds \int_{\Pi_v^-\cap \Om } G_j(|x|, u_1(x), |\nabla u_1(x)|) \, dx .
$$
As previously, it is not hard to see that $\psi _j$ is a continuous odd mapping 
on $S^{N-1}$, $ 1 \leq j \leq k$. 
In particular, the restriction of $ \Psi = (\psi _1, \dots, \psi _k)$ to 
$S^{N-1} \cap \Pi_{e_1}$ is a continuous odd mapping from this space to $\R^k$. 
Since  $S^{N-1} \cap \Pi_{e_1}$ can be identified to $S^{N-2}$ and 
$ k \leq N-2$, 
we may use the Borsuk-Ulam theorem again and we infer that there exists 
$ e_2 \in S^{N-1} \cap \Pi_{e_1} $ such that $ \Psi ( e_2) = 0$, 
i.e. $ \Pi_{e_2} $ splits the constraints in two for the minimizer $ u_1$. 
We denote $ u_{1,1} = (u_1) _{\Pi_{e_2} ^-} $ and 
$ u_{1,2} = (u_1) _{\Pi_{e_2} ^+} $ the  functions obtained from $u_1$ by the reflection procedure (1). 
Arguing as previously, we infer that $ u_{1,1} $ and $ u_{1,2}$ 
belong to $\Xo $, satisfy the constraints and are minimizers for problem 
($\Po$).
Moreover, they are symmetric with respect to $ \Pi_{e_1}$ and with respect to 
$ \Pi_{e_2}$.
Next we use the following: 

\bigskip

\noindent
{\bf Lemma 3. } \it Let $w\in \Xo $  be a minimizer for ($\Po$).
Assume that {\bf A1, A2 } are satisfied and there exists a vector 
subspace $V$ of $\R^N$ of dimension  $m \leq N-2$ such that 
 any hyperplane containing $V$ splits the constraints in two for $w$. 
Then $w$ is radially symmetric with respect to $V$. 

\rm 

\bigskip

{\it Proof. } 
Let $\Bo _1 =\{ b_1, \dots, b_{m} \}$ be an orthonormal basis 
in $V$. 
Fix a hyperplane $\Pi $ containing $V$. 
We extend $\Bo _1  $ to an orthonormal basis $ \Bo = \{b_1, \dots, b_{N} \}$ 
in $\R^N$ in such a way that $ \Pi = \Pi_{b_N} = b_N ^{\perp}$. 
We denote by $ (x_1, \dots, x_N)$ the
coordinates of a point $ x $ with respect to $ \Bo$.
Let  $w_1 = w_{\Pi_{b_N}^-} $ and $ w_2 = w_{\Pi_{b_N}^+} $.
Clearly $ w_1, w_2 \in \Xo $ by {\bf A1}. 
By the assumption of Lemma 3, $ \Pi_{b_N} $ splits the constraints in two for 
$ w$ and this implies that $w_1$ and $w_2$ satisfy the constraints.
As before we have $ E(w_1) \geq E(w)$, $ E(w_2) \geq E(w)$
and $ E(w_1 )+ E(w_2) = 2 E(w)$, 
thus necessarily $ E(w_1) = E(w_2) = E(w)$ and $ w_1$, $w_2$ are also minimizers. By {\bf A2} we have $ w, w_1, w_2 \in C^1(\Om)$. 
Since $ w_1$ is  symmetric with respect to the $ x_N$ variable, we have 
$\frac {\p w_1 }{\p x_N} ( x_1, \dots, x_{N-1}, 0) = 0 $ 
  whenever $ ( x_1, \dots, x_{N-1}, 0)  \in \Om$.
But $ w(x) = w_1 (x) $ for $ x_N <0$,  therefore
$$
\begin{array}{l}
\ds \frac {\p w }{\p x_N} ( x_1, \dots, x_{N-1}, 0) = 
\ds \lim_{s \uparrow 0 } \frac {\p w }{\p x_N} ( x_1, \dots, x_{N-1}, s) 
\\
\\
= 
\ds \lim_{s \uparrow 0 }  \frac {\p w_1 }{\p x_N} ( x_1, \dots, x_{N-1}, s) = 
\frac {\p w_1 }{\p x_N} ( x_1, \dots, x_{N-1}, 0) = 0
\end{array}
\leqno{(3)}
$$
for  $ ( x_1, \dots, x_{N-1}, 0) \in \Om$, 
i.e. the  derivative of $w$ in  the direction orthogonal
 to $ \Pi $ vanishes on $ \Om \cap \Pi$. 
Thus we have proved that for any hyperplane $ \Pi $  containing $V$, we have 
$$
\frac{\p w}{\p n } (x) = 0
\qquad \mbox{ for any } x \in \Om \cap \Pi, 
\mbox{ where $ n $ is the unit normal to } \Pi.
\leqno{(4)}
$$

We pass to spherical coordinates in the last $ N-m$ variables in $\R^N$, i.e. we use variables $ (r, \theta _1, \dots, \theta _{N-m -1})$ 
instead of $ (x_{m+1}, \dots, x_N)$, 
where $ r = \left( x_{N -m +1} ^2 + \dots + x_N ^2 \right)^{\frac 12} $
and $\theta _1, \dots \theta _{N-m -1}$ are the angular variables. 
Then (4) is equivalent to $ \frac{ \p w}{\p \theta _j} = 0$ on $ \Om $ for 
$ j =1, \dots, N-m -1$. 
We infer that $ w$ does not depend on $ \theta _1, \dots, \theta _{N-m +1}$,
 i.e. there exists some function $ \tilde{w} $ depending only on 
$ x_1, \dots, x_m, r$ such that 
$ w( x_1, \dots, x_N) = \tilde{w}(x_1, \dots, x_m, r)$ on $ \Om $ and Lemma 3 is proved. 
\hfill
$\Box$

\bigskip

Now  come back to the proof of Theorem 1. 
Clearly, any  $ x \in \R^N$ has a unique decomposition 
 $ x = x_1 e_1 + x_2 e_2 + x'$, where 
$ x_1, x_2 \in \R$ and $ x' \in \{ e _1, e_2 \}^{\perp }$. 
Since $ u_{1,1 }$ and $ u_{1,2}$ are symmetric with respect to 
$ \Pi_{e_1}$ and with respect to 
$ \Pi_{e_2}$, we have 
$ u_{1, i} ( x_1 e_1 + x_2 e_2 + x') = u_{1,i}(x_1 e_1 - x_2 e_2 + x') 
=  u_{1, i} ( -x_1 e_1 - x_2 e_2 + x')$. 
Let $ \Pi $ be a hyperplane containing $\{ e_1, e_2 \} ^{\perp }$. 
It is obvious that the transform 
$ x_1 e_1 + x_2 e_2 + x' \longmapsto  -x_1 e_1 - x_2 e_2 + x' $ 
is a one-to-one correspondence between $ \Pi ^+ $ and $ \Pi ^-$ 
and a simple change of variables gives
$$
\ds \int_{\Pi ^+ \cap \Om } G_j(|x|, u_{1, i} (x), |\nabla u_{1,i}(x)|) \, dx 
= \ds \int_{\Pi ^-\cap \Om } G_j(|x|, u_{1,i}(x), |\nabla u_{1,i} (x)|) \, dx , 
\quad
 j = 1, \dots, k, 
$$
hence $ \Pi $ splits the constraints in two for $ u_{1,i}$, $ i =1,2$. 
By Lemma 3, we infer that $ u_{1,i}$
 are radially symmetric with respect to 
$\{ e_1, e_2 \} ^{\perp }$, i.e.  
$ u_{1,i}( x_1 e_1 + x_2 e_2 + x' ) = \tilde{ u}_{1,i} 
( \sqrt{ x_1 ^2 + x_2 ^2} , x') $ for some functions $ \tilde{ u}_{1,1} $ and
$ \tilde{ u}_{1,2} $.
On the other hand, $ u_{1,1}(x) = u_1 (x) = u_{1,2}(x) $ for any 
$ x \in \Pi _{e_2 } \cap \Om $,
that is $\tilde{ u}_{1,1} ( |x _1| , x') =  \tilde{ u}_{1,2} ( |x_1| , x')$
whenever $ x_1 e _1 + x' \in \Om $ . 
We conclude that necessarily 
$\tilde{ u}_{1,1}  =  \tilde{ u}_{1,2} $
and $ u_{1,1} ( x) = u_1(x)  = u_{1,2}(x) $ for any $ x \in \Om $, 
thus $ u_1 $ is radially symmetric with respect to 
$\{ e_1, e_2 \} ^{\perp }$. 

\medskip

Similarly there exists $ v_ 2 \in S^{N-1} \cap e_1^{\perp } $ such that 
$ \Pi _{v_2 }$ splits the constraints in two for $ u_2 $ and we infer that 
$u_2 $ is radially symmetric with respect to 
$\{ e_1, v_2 \} ^{\perp }$. 
We use this information together with the fact that $ u_1 = u = u_2$ on 
$ \Om \cap \Pi _{e_1 }$ 
to prove the symmetry of $u$.

If $ v_2 $ is colinear to $ e_2$, i.e. $ v_2 = \pm e_2$,  we may assume that 
$ v_2 = e_2$. 
Using the symmetry of $ u_1, $ $ u_2 $ and the fact that 
$ u_1 = u = u_2 $ on $ \Om \cap \Pi_{e _1}$, we obtain as above that 
 $u_1 = u_2 = u $ on $\Om$, hence $u$ is radially symmetric with respect to 
$\{ e_1, e_2 \} ^{\perp }$.
 
If $ v_2 $ and $ e_2 $ are not colinear,  $\mbox{Span}\{ e_1, e_2, v_2 \}$ 
is a three-dimensional subspace. Let $ \{ e_4, \dots , e_N \}$ be an orthonormal
basis in $ \{ e_1, e_2, v_2 \} ^{\perp} $. 
We choose $ e_3 $ and $ v_3 $ in such a way that 
$ \Bo = \{ e_1, e_2, e_3, \dots, e_ N \}$ and 
$ \Bo '  = \{ e_1, v_2, v_3, e_4,  \dots, e_ N \}$  
are orthonormal basis in $ \R^N$ with the same orientation.
Then there exists $ \te \in (0, \pi ) \cup ( \pi , 2 \pi ) $ such that 
$ v_2 = \cos \te\,  e _2 + \sin \te \, e_3$ and $ v_3 = - \sin \te \, e_2 + \cos \te \, e_3$. 
Given a point $ x \in \R^N$, we denote by $ ( x_1, x_2, \dots, x_N )$ its coordinates with respect to $ \Bo$.  
It is clear that $(x_1, y_2 = \cos \te \, x_2 + \sin \te \, x_3, 
y_3 = - \sin \te \, x_2 + \cos \te \, x_3 , x_4, \dots, x_N)$ 
are the coordinates of $ x$ with respect to $ \Bo '$.

Fix $ r e_3 +  \sum_{j = 4 }^N x_j e_j \in \Om \cap e_1^{\perp} $ 
and denote 
$$
\ph (t) = \ph _{r, x_4, \dots, x_N} (t) = 
u( r \cos t \, e _2 + r \sin t \, e _3 + \ds \sum_{j = 4 }^N x_j e_j ).
$$
Clearly, $ \ph $ is $ C^1 $ and $ 2\pi -$periodic on $ \R$. 
Since the restriction of 
$ u = u_1 $  to $\Om \cap e_1^{\perp} $ 
is symmetric with respect to $ \R e_2$,  we get 
$$
\ph (t) = 
u( - r \cos t \, e _2 + r \sin t \, e _3 + \ds \sum_{j = 4 }^N x_j e_j )
= \ph ( \pi -t). 
\leqno{(5)}
$$
The restriction of  $ u = u_2 $ to $\Om \cap e_1^{\perp} $
is also symmetric with respect to $ \R v_2$, therefore 
$$
\begin{array}{l}
\ph (t) = u ( r ( \cos t \cos \te + \sin t \sin \te ) \, v_2 
+ r(\sin t \cos \te - \cos t \sin \te ) v_3 +  \sum_{j = 4 }^N x_j e_j)
\\
\\
= u_2 ( r \cos( t - \te ) \, v_2 + r \sin ( t - \te ) \, v_3 
+ \sum_{j = 4 }^N x_j e_j ) 
\\
\\
= u_2 (-r \cos( t - \te ) \, v_2 + r \sin ( t - \te ) \, v_3 
+ \sum_{j = 4 }^N x_j e_j ) 
\\
\\
=u_2( r \cos( \pi - ( t - \te )) \, v_2 + r \sin( \pi - ( t - \te ))\, v_3 
+ \sum_{j = 4 }^N x_j e_j)
\\
\\ 
= \ph ( \pi + 2 \te  - t ) = \ph ( t - 2 \te ) 
\qquad \mbox{ by (5)}. 
\end{array}
\leqno{(6)}
$$
Hence any of the functions $ \ph _{r, x_4, \dots, x_N}$ admits  
$ 2 \pi $ and $ 2 \te $ as periods. 
The following situations may occur: 

\medskip

Case 1: 
$ \frac{ \te }{\pi } \in \R \setminus \Q$. 
The set $ \{ 2n \te + 2 k \pi \; | \; n, k \in \Z \} $
is dense in $ \R$ and any number in this set is a period for  
$ \ph _{r, x_4, \dots, x_N}$. 
Since $ \ph _{r, x_4, \dots, x_N}$ is continuous, we infer that it is constant. This is equivalent to 
$ u (\sum_{j = 2 }^N x_j e_j )
= u( \sqrt{ x_2 ^2 + x_3 ^2 } \, e_2 + \sum_{j = 4 }^N x_j e_j)$  whenever 
$ \sum_{j = 2 }^N x_j e_j \in \Om \cap e_1^{\perp}$.
With the above notation, 
using the symmetry properties of $ u_1 $ and $ u_2 $ 
we have for any $ x \in \Om$, 
$$
u_1 (x)= u_1(  \sqrt{ x_1 ^2 + x_2 ^2 } \, e_2 + \sum_{j = 3 }^N x_j e_j)
= u (  \sqrt{ x_1 ^2 + x_2 ^2 + x_3 ^2 } \, e_2 + \sum_{j = 4 }^N x_j e_j)
$$
and
$$
\begin{array}{l}
u_2 (x)= u_2(  \sqrt{ x_1 ^2 + y_2 ^2 } \, v_2 +  y_3 v_3 +
\sum_{j = 4 }^N x_j e_j)
= u(  \sqrt{ x_1 ^2 + y_2 ^2 } \, v_2 +  y_3 v_3 +
\sum_{j = 4 }^N x_j e_j)
\\
\\
= u (  \sqrt{ x_1 ^2 + y_2 ^2 + y_3 ^2 } \, e_2 + \sum_{j = 4 }^N x_j e_j)
= u (  \sqrt{ x_1 ^2 + x_2 ^2 + x_3 ^2 } \, e_2 + \sum_{j = 4 }^N x_j e_j).
\end{array}
$$
Consequently $ u = u_1 = u_2 $ on $ \Om $ and $ u$ is radially symmetric with 
respect to $ \{e_1, e_2, e_3 \}^{\perp}$. 

\medskip

Case 2: $ \frac{ \te }{\pi} = \frac kn$ where $ k, n $ are relatively prime integers, $ k$ is odd and $n$ is even, say $ k = 2k_1 + 1 $ and $ n = 2 n_1$. 
Then $ \pi =2 n_1 \te - 2 k_1 \pi $ is also a period for 
$ \ph_{r, x_4, \dots, x_N}$ and this implies
$$
 u( \sum_{j = 2 }^N x_j e_j) 
= u( - x_2 e_2 - x_3 e_3 + \sum_{j = 4 }^N x_j e_j) 
\qquad
\mbox{ whenever } \sum_{j = 2 }^N x_j e_j \in \Om \cap e_1 ^{\perp}.
\leqno{(7)}
$$

From the symmetry of $ u_1$  and (7) we get for $ x_1 \leq 0$:
$$
\begin{array}{l}
u (  \sum_{j = 1}^N x_j e_j)
 = u( \sqrt{ x_1 ^2 + x_2 ^2}\,  e_2 +  \sum_{j = 3 }^N x_j e_j) 
\\
\\
= u(  - \sqrt{ x_1 ^2 + x_2 ^2}\,  e_2 - x_3 e_3 +  \sum_{j = 4 }^N x_j e_j)  
= u( x_1 e_1 - x_2 e_2 - x_3 e_3 + \sum_{j = 4 }^N x_j e_j)  .
\end{array}
\leqno{(8)}
$$
Using the symmetry of $ u_2$ and (7), we infer 
 that (8) also  holds  for $ x_1 \geq 0$.   
Let $\Pi $ be a hyperplane containing $ \{ e_1, e_4, , \dots , e_N \}$. 
It is clear that the mapping 
$ \sum_{j = 1}^N x_j e_j \longmapsto  
x_1 e_1 - x_2 e_2 - x_3 e_3 + \sum_{j = 4 }^N x_j e_j $ 
is a linear isometry between $ \Pi ^+ $ and $ \Pi ^-$. 
Then (8) and a simple change of variables show that 
$$
\int_{\Pi ^+ \cap \Om } G_{\ell} ( |x|, u(x), |\nabla u(x)|) \; dx 
= \int_{\Pi ^- \cap \Om } G_{\ell} ( |x|, u(x), |\nabla u(x)|) \; dx , 
$$
for $ \ell = 1, \dots, k$, 
i.e. $ \Pi $ splits the constraints in two for $ u$. 
Since $u$ is a minimizer, by Lemma 3 we infer that $u$ is radially symmetric with respect to $\mbox{ Span}\{ e_1, e_4, \dots, e_N \}$. 
In particular, the restriction of $u$ to $ \Om \cap e_1^{\perp }$ is radially symmetric with respect to $\mbox{ Span}\{  e_4, \dots, e_N \}$. 
As in case 1, this implies that $u$ is radially symmetric with respect to 
$\mbox{ Span}\{  e_4, \dots, e_N \}$. 

\medskip

Case 3: $ \frac{ \te }{\pi} = \frac kn$ where $ k, n $ are relatively prime integers, $ k$ is even  and $n$ is odd, say $ k = 2k_1  $ and $ n = 2 n_1+1$. 
Then $ \te =  2 k_1 \pi - 2 n_1 \te$ is  a period for 
$ \ph_{r, x_4, \dots, x_N}$.
By (5) we get 
$ \ph_{r, x_4, \dots, x_N} (t) = \ph_{r, x_4, \dots, x_N} ( \pi - t) = 
\ph_{r, x_4, \dots, x_N}( \te + \pi -t)$. This means that 
for $ \sum_{j = 2}^N x_j e_j \in \Om $ we have
$$
u (\sum_{j = 2}^N x_j e_j   ) 
= u ( - (x_2 \cos \te  + x_3 \sin \te ) e_2 + 
( - x_2 \sin \te + x_3 \cos \te ) e_3 
+  \sum_{j = 4}^N x_j e_j ) .
\leqno{(9)}
$$
In other words, for fixed $ x'' \in \mbox{ Span}\{ e_4, \dots, e_N \}$, 
the function $ x_2 e_2 + x_3 e_3 \longmapsto u ( x_2 e_2 + x_3 e_3 + x'') $
is symmetric with respect to $ \R w$, 
 where 
$ w = \cos( \frac{ \te + \pi}{2} ) e_2 + \sin ( \frac{ \te + \pi}{2} ) e_3$. 
Note that the symmetry of $\mbox{ Span}\{e_1, e_2, e_3 \}$ 
with respect to $ \R w$ is a linear isometry of matrix
$ A = \left( \begin{array}{ccc}
-1 & 0 & 0 \\
0 & - \cos \te  & - \sin \te \\
0 & - \sin \te  & \cos \te 
\end{array}
\right)
$
with respect to the basis $ \{ e_1, e_2, e_3 \}$. 
We show that for any $ x \in \Om $ we have
$$
u(x) = u( Sx),
\leqno{(10)}
$$
where 
$ 
Sx = - x_1 e_1 - (x_2 \cos \te  + x_3 \sin \te ) e_2 + 
( - x_2 \sin \te + x_3 \cos \te ) e_3 
+  \sum_{j = 4}^N x_j e_j.
$
It suffices to consider the case $ x_1 \leq 0$. 
By using the symmetry of $u_1$, $u_2$ and (9) we get 
$$
\begin{array}{l}
u(x) = u_1 (x) = u( \sqrt{ x_1 ^2 + x_2 ^2}\, e_2 + \sum_{j = 3}^N x_j e_j ) 
\\
\\
= u ( - (\sqrt{ x_1 ^2 + x_2 ^2} \cos \te + x_3 \sin \te ) e_2 
+ ( - \sqrt{ x_1 ^2 + x_2 ^2} \sin \te + x_3 \cos \te ) e_3   
+ \sum_{j = 4}^N x_j e_j )
\end{array}
$$
and 
$$
\begin{array}{l} 
u( Sx ) = u_2 ( Sx ) = u_2( - x_1 e_1 - x_2 v_2 + x_3 v_3 
+ \sum_{j = 4}^N x_j e_j )
\\
\\
= u_2( - \sqrt{ x_1 ^2 + x_2 ^2}\, v_2 +  x_3 v_3 
+ \sum_{j = 4}^N x_j e_j )
\\
\\ 
= u (- \sqrt{ x_1 ^2 + x_2 ^2} ( \cos \te \, e_2 + \sin \te \, e_3)  
+ x_3 ( - \sin \te \, e_2 + \cos \te \, _3 ) +  \sum_{j = 4}^N x_j e_j ), 
\end{array}
$$
hence $ u(x) = u ( Sx)$. 
Let $ \Pi $ be a vector hyperplane containing $w, e_4, \dots, e_N$. 
It is easy to see that $ S$ is a linear isometry of $ \R^N$ mapping 
$ \Om \cap \Pi ^- $ onto $ \Om \cap \Pi ^+$. 
Using (10) and a change of variables, 
we find that $\Pi$ splits the constraints in two for $u$. 
By Lemma 3 we infer that $u$ is radially symmetric with respect to 
$ \mbox{ Span} \{ w, e_4, \dots, e_N \}$.

In fact, since $u_1$ is radially symmetric with respect to 
$ \mbox{ Span} \{ e_3, e_4, \dots, e_N \}$ and  
$ \mbox{Span} \{ w, e_4, \dots, e_N \}$, 
it can be proved that $u_1 $ is radially symmetric with respect to 
$ \mbox{Span} \{  e_4, \dots, e_N \}$.
Similarly $ u_2 $ is radially symmetric with respect to 
$ \mbox{Span} \{  e_4, \dots, e_N \}$ 
and then it is clear that $u$ has the same property. 
We omit the proof because we will not  make use of this observation.

\medskip

Case 4: $ \frac{ \te }{\pi} = \frac kn$ where $ k, n $ are relatively prime odd integers,  say $ k = 2k_1 + 1 $ and $ n = 2 n_1+1$. 
Then $ \te - \pi =  2 k_1 \pi - 2 n_1 \te$ is  a period for 
$ \ph_{r, x_4, \dots, x_N}$.
By (5) we have 
$ \ph_{r, x_4, \dots, x_N} (t) = \ph_{r, x_4, \dots, x_N} ( \pi - t) = 
\ph_{r, x_4, \dots, x_N}( \te  -t)$, that is 
$$
u(x) = 
u( (x_2 \cos \te + x_3 \sin \te ) e_2 + ( x_2 \sin \te - x_3 \cos \te ) e_3  
+ \sum_{j = 4}^N x_j e_j )
\leqno{(11)}
$$
for any $ x = \sum_{j = 2}^N x_j e_j  \in \Om \cap e_1^{\perp }$. 
Proceeding as in case 3, we prove that $u$ is radially symmetric with 
respect to $ \mbox{ Span} \{ w' , e_4, \dots, e_N \}$, 
where $ w' = \cos \frac{ \te }{2} \, e_2 + \sin \frac{ \te }{2} \, e _3$. 
(In fact, it can be proved that $u$ is radially symmetric with 
respect to $ \mbox{ Span} \{ e_4, \dots, e_N \}$).

\medskip

Note that in either case it follows that $u$ is symmetric with respect to 
$ \Pi _{e_1}$. 
Thus we have proved that whenever $ e _1 \in S^{N-1}$ 
 satisfies $ \Phi ( e_1 ) = 0$, 
$u$ is symmetric with respect to $ \Pi _{e_1}$. 
Assume that $ e_1, \dots , e_{\ell} \in S^{N-1 }$ are mutually orthogonal,  satisfy $ \Phi ( e_1) = \dots = \Phi ( e_{\ell} ) = 0 $ and 
$ \ell \leq N - k -1$. 
It is clear that $ S_{\ell} = S^{N-1} \cap \{ e_1, \dots, e_{\ell} \}^{\perp }$ 
can be identified to $ S^{ N-\ell -1 }$ and the restriction of 
$ \Phi $ to $ S_{\ell } $ is an odd, continuous function from $ S_{\ell } $ to 
$ \R^k$. 
Using the Borsuk-Ulam theorem we infer that there exists 
$ e_{\ell + 1 } \in S_{\ell } $ such that $ \Phi ( e_{\ell +1} ) = 0$. 
By induction it follows that there exist $ N-k$ mutually orthogonal vectors 
$ e_1, \dots , e_{N -k } \in S^{N-1 }$ such that 
$ \Phi ( e_1) = \dots = \Phi ( e_{N -k} ) = 0 $. 
We complete this set to an orthonormal basis $ \{ e_1, \dots, e_N \}$
in $ \R^N$. 
We already know that $u$ is symmetric with respect to 
any of the hyperplanes $ \Pi _{e_1}, \dots, \Pi_{e_{N-k}}$. 
In particular, for $ x = \sum_{j = 1}^N x_j e_j \in \Om $ we have
$$
u(x) = u( - x_1 e_1 + \sum_{j = 2}^N x_j e_j )
= \dots 
= u( -  \sum_{j = 1}^{N-k} x_j e_j + \sum_{j = N -k +1}^N x_j e_j).
\leqno{(12)}
$$  
Let $ \Pi $ be a (vector)  hyperplane containing $ e_{N-k +1}, \dots, e_N$. 
It is clear that the mapping 
$ \sum_{j = 1}^N x_j e_j \longmapsto 
-  \sum_{j = 1}^{N-k} x_j e_j + \sum_{j = N -k +1}^N x_j e_j $ 
is a linear isometry between $ \Pi ^+ $ and $ \Pi ^-$. 
Using (12), we infer that $ \Pi $ splits the constraints in two for $ u$.
By Lemma 3, $u$ is radially symetric with respect to 
$ \mbox{ Span}\{ e_{N-k +1}, \dots, e_N \}$.

\medskip

The  case $k=0$ is much simpler. 
Problem ($\Po $) consists in   minimizing $E$ on $\Xo $ without constraints. 
Assume that $ u$ is a minimizer. 
Let $ \Pi $ be a hyperplane containing the origin and let $ u_{\Pi ^-}$, 
$ u_{\Pi ^+ }$ be the two functions obtained from $u$ as in (1). 
By {\bf A1 } we have $ u_{\Pi ^-}, u_{\Pi ^+} \in \Xo $, thus 
$ E( u_{\Pi ^-}) \geq  E( u) $ and $ E( u_{\Pi ^+}) \geq  E( u) $. 
On the other hand,  $ E( u_{\Pi ^-}) + E( u_{\Pi ^+}) = 2E(u)$, thus necessarily $ E( u_{\Pi ^-}) = E( u_{\Pi ^+}) = E(u)$ 
and $ u_{\Pi ^-}$, 
$ u_{\Pi ^+ }$ are also minimizers. 
As in the proof of Lemma 3, this implies $ \frac { \p u }{\p n } ( x) = 0$ 
for any $ x \in \Om \cap \Pi$, where $ n $ is the unit normal to $ \Pi$. 
Then passing to spherical coordinates, as in Lemma 3, we see that $u$ does not depend on the angular variables, i.e. $u$ is a radial function. 
\hfill $\Box$

\bigskip

{\it Proof of Theorem 2. } 
For $ v \in S^{ N-1} $ and $ t \in \R$ we denote by $ \Pi_{v,t } $ the affine hyperplane $\{ x \in \R^N \; | \; ( x - tv ) . v = 0 \} $ and by 
$ \Pi_{v,t } ^ + = \{ x \in \R^N \; | \; ( x - tv ) . v > 0 \}$, respectively 
$ \Pi_{v,t } ^ - = \{ x \in \R^N \; | \; ( x - tv ) . v < 0 \}$
the two half-spaces determined by $ \Pi_{v,t } $. 
It is clear that $\Pi_{-v,-t }^ -  = \Pi_{v,t }^+ $.
For $ j =1, \dots, k$, 
 we define $ \tilde{ \psi } _j :S^{N-1} \times \R \lra \R$ 
by 
$$
\tilde{ \psi } _j (v,t) =
\ds \int_{\Pi_{v ,t }^+  } G_j( u(x), |\nabla u(x)|) \, dx 
- \ds \int_{\Pi_ {v, t}^-  } G_j( u(x), |\nabla u(x)|) \, dx .
$$
 Since $ G_j ( u, |\nabla u |) \in L^1(\R^N)$, 
it is a simple consequence of Lebesgue's dominated convergence theorem that 
$ \tilde{\psi }_j$ is continuous on $ S^{N-1} \times \R $.
It is obvious that $ \tilde{\psi }_j ( -v , -t ) = - \tilde{\psi }_j( v, t)$.

We claim that $ \ds \lim_{t \ra \infty } \tilde{\psi }_j( v, t) 
= - {\ds \int_{\R^N} } G_j( u(x), |\nabla u (x)|) \, dx = - \la _j $ uniformly
with respect to $ v \in S^{N-1}$.
Indeed, fix $ \e > 0$. 
There exists $ R > 0$ such that 
$$
{\ds \int_{\R^N  \setminus B(0,R) } } | G_j( u(x), |\nabla u (x)|) |  \, dx  
< \frac{ \e}{2}.
$$
For any $ v \in S^{N-1}$ and $ t  > R$ we have 
$\Pi_{v,t } ^ +  \subset \R^N  \setminus B(0,R) $, therefore
$$
\bigg\vert 
\tilde{\psi }_j( v, t)  + 
{\ds \int_{\R^N} } G_j( u(x), |\nabla u (x)|) \, dx 
\bigg\vert 
= 2 \bigg\vert  
{\ds \int_{\Pi_{v,t } ^ + }}G_j( u(x), |\nabla u (x)|) \, dx  
\bigg\vert 
< \e
$$
and the claim is proved.
It is clear that 
$ \ds \lim_{t \ra - \infty } \tilde{\psi }_j( v, t) = \la _j $ 
uniformly in $ v \in  S^{N-1}$.

We denote $P= ( 0, \dots, 0 , 1) \in \R^{N +1} $, 
$S = ( 0, \dots, 0 , -1) \in \R^{N +1} $ 
and we define $ \psi _j : S^N \lra \R$ by 
$$
\psi _j ( x_1, \dots, x_N, x_{N+1}) = 
\tilde{\psi}_j
\left( \frac{( x_1, \dots, x_N )}{|( x_1, \dots, x_N )| } , 
\frac{x_{N+1}}{1 - |x_{N+1}|} \right) 
$$
if $ ( x_1, \dots, x_N, x_{N+1}) \not\in \{ P, S \}$, respectively 
$ \psi _j ( P) = - \la _j$ and $ \psi _j (S) = \la _j$. 
Then $ \psi _j$ is an odd, continuous function on $ S^N$. 

\medskip

Consider first the case $ 1 \leq k \leq N - 2$. 
It follows from Theorem 1 that there exist two orthogonal vector subspaces
$ V_1 $ and $ V_2 $ such that $\mbox{ dim}(V_1) = k$, $ V_1 \oplus V_2 = \R^N$ and $u$ is radially symmetric with respect to $ V_1$. 
The set $ {\mathbf{S}} = \{ ( y_1, \dots, y_N , y_{N+1 } ) \in S^N \; | \; 
( y_1 , \dots, y_N) \in V_1  \} $ 
can be identified to $ S^k$. 
Since the restriction of $ \Psi = (\psi _1, \dots, \psi _k) $ to 
$ {\mathbf{S} } \simeq S^k $ is continuous, odd, $\R^k -$valued, by the Borsuk-Ulam 
theorem we infer that there exists 
$ y^* = ( y_1 ^*, \dots, y_N^*, y_{N+1}^* )\in {\mathbf{S}}$ such that 
$ \psi ( y^*) = 0$. 
We cannot have $ y^* = S $ or $ y^* = P$ because 
$ \psi (S) = - \psi(P) = ( \la _1, \dots, \la _N) \neq 0$. 
Denote $ e_k = \frac{( y_1^* , \dots, y_N ^*)}{|(  y_1 ^* , \dots, y_N ^*)|}$ 
and $ t = \frac{y_{N+1}^*}{1 - | y_{N+1}^*|}$. 
Then $ e_k \in V_1$, $|e_k| = 1$ and $ \tilde{\psi} _j( e_k, t) = 0$ for 
$ j =1, \dots, k$, i.e. $\Pi_{e_k, t }$ splits the constraints in two for $u$. 
Choose $e_i $, $ i =1, \dots, N$, $ i \neq k$ in such a way that 
$ \{ e_1, \dots, e_{k-1}, e_k \}$ and 
 $ \{ e_{k+1}, \dots, e_N \}$ are orthonormal basis in $ V_1$, respectively 
in $ V_2$. 
Denote $ u_*(x) = u ( x - t e_k)$. 
It is clear that $ u_*$ is a minimizer for ($\Po '$), it is radially symmetric with respect to $V_1 $ and the hyperplane $ e_k^{\perp } = \Pi_{e_k, 0} $ 
splits the constraints in two for $ u_*$. Arguing exactly as in the proof of Theorem 1, we see that $u_* $ is symmetric with respect to $  e_k^{\perp }  $. 
Using this fact and the radial symmetry with respect to $V_1$, we get
$$
u_*(\sum_{i =1}^N x_i e_i )  
= u_* ( \sum_{i =1}^k x_i e_i - \sum_{i =k+ 1}^N x_i e_i) 
= u_* ( \sum_{i =1}^{k-1} x_i e_i - \sum_{i =k }^N x_i e_i).
\leqno{(13)}
$$
By (13) we infer that any (vector) hyperplane containing 
$ e_1, \dots, e_{k-1}$  splits the constraints in two for $u_*$. 
Then Lemma 3 implies that $u_*$ is radially symmetric with respect to 
$\mbox{ Span}\{ e_1, \dots, e_{k-1} \}$,
consequently $ u$ is radially symmetric with respect to the affine subspace 
$ t e_k + \mbox{Span}\{ e_1, \dots, e_{k-1} \}$.

\medskip

Now consider the case $ k = N-1$. 
As above, there exists $ y^* = ( y_1 ^*, \dots, y_N^*, y_{N+1}^* )
\in S^N \setminus \{ S, P \}$ such that $ \psi ( y^*) = 0$. 
Denoting $ e_1 = \frac{( y_1 ^*, \dots, y_N ^*)}{|(  y_1 ^*, \dots, y_N ^*)|}$ 
and $ t_1 = \frac{y_{N+1}^*}{1 - | y_{N+1} ^*|}$, this means that 
$ \Pi_{e_1, t_1 }$ splits the constraints in two for $u$. 
Let $ u_1 = u_{\Pi_{e_1, t_1}^-} $ and
$ u_2 = u_{\Pi_{e_1, t_1}^ +} $.
It is clear that $u_1$, $u_2 $  are also minimizers for ($\Po '$).
Since $\{ (y_1, \dots, y_{N+1}) \in S^N \; | \; 
( y_1, \dots , y_N) \perp e_1 \}$ is 
  homeomorphic to $ S^{N-1 }$ and there are exactly $ N-1$ constraints,  it is 
 possible to restart the 
prevoius process with $ u_1 $ instead of $u$. 
We infer that there exists $e_2 \in e_1^{\perp}$, $|e_2 | =1$ and $ t_2 \in \R$
such that $ \Pi_{e_2, t_2 }$ splits the constraints in two for $ u_1$. 
Putting $ u_{1,1} = (u_1)_{\Pi_{e_2, t_2}^-} $ and
$ u_{1,2} = (u_1) _{\Pi_{e_2, t_2}^ + } $, we see that $ u_{1,1} $ and 
$u_{1,2} $ are minimizers for ($\Po '$)
and are symmetric  with respect to $\Pi_{e_1, t_1 }$ and 
$\Pi_{e_2, t_2 }$. 
It follows that $ \tilde{u}_{1,1} = u_{1,1} ( \cdot - t_1 e_1 - t_2 e_2 ) $ 
and $ \tilde{u}_{1,2} = u_{1,2} ( \cdot - t_1 e_1 - t_2 e_2 ) $ 
minimize ($\Po '$) and are symmetric with respect to $ e_1^{\perp } $ and 
$ e_2^{\perp }$. 
Therefore any (vector) hyperplane in $ \R^N$ containing 
$ \{ e_1, e_2 \}^{\perp }$ 
splits the constraints in two for $ \tilde{u}_{1,1}$ and for $ \tilde{u}_{1,2}$
and using Lemma 3 we infer that  $ \tilde{u}_{1,1}$ and  $ \tilde{u}_{1,2}$
are radially symmetric with respect to $ \{ e_1, e_2 \}^{\perp }$.
Since $ \tilde{u}_{1,1}= \tilde{u}_{1,2}$  on $ \Pi_{e_2, 0} =e_2 ^{\perp }$, we have 
necessarily $ \tilde{u}_{1,1}= \tilde{u}_{1,2}$  on $\R^N$. 
Therefore $ u_1 = \tilde{u}_{1,1} ( \cdot + t_1 e_1 + t_2 e_2 ) $ 
is radially symmetric with respect to the affine subspace
$ t_1 e_1 + t_2 e_2 + \{ e_1, e_2 \}^{\perp }$.

Similarly we prove that there exist $ v_2 \in e_1^{\perp}$, $|v_2 | = 1$ and 
$ s_2 \in \R$ 
such that $ u_2 $ is radially symmetric with respect to the affine subspace
$ t_1 e_1 + s_2 v_2 + \{ e_1, v_2 \}^{\perp }$.
Of course, nothing guarantees {\it \`a priori } that 
$( e_2, t_2) = \pm ( v_2, s_2)$. The following situations may occur:

\medskip

Case 1: $ e_2$ and $ v_2 $ are colinear. Then we may assume that $ e_2 = v_2$. 
There are two subcases: 

a) $ t_2 = s_2$. Then $ u_1 ( \cdot - t_1 e_1 - t_2 e_2 ) $ and 
$ u_2 ( \cdot - t_1 e_1 - t_2 e_2 ) $ are both radially symetric with respect to $\{ e_1, e_2 \}^{\perp } $ and are equal on $ e_1^{\perp }$. 
We conclude that $  u_1 ( \cdot - t_1 e_1 - t_2 e_2 ) =
u_2 ( \cdot - t_1 e_1 - t_2 e_2 ) $, 
thus $ u = u_1 = u_2 $ is radially symmetric with respect to 
$  t_1 e_1  + t_2 e_2 + \{ e_1, e_2 \}^{\perp } $.

b)  $ t_2 \neq s_2$, say $ s_2 > t_2$. The symmetry of $ u_1 $ and $ u_2 $ imply that there exist some functions $ \tilde{u}_1$, $ \tilde{u}_2$
defined on $ [0, \infty ) \times \{ e_1, e_2 \}^{\perp } $ such that 
$$
\begin{array}{l}
u_1(x_1 e_1 + x_2 e_2 + x') = 
\tilde{u}_1 (\sqrt{ ( x_1 - t_1 )^2 + ( x_2 - t_2 )^2} \, , x' ) 
\\
\\
u_2(x_1 e_1 + x_2 e_2 + x') = 
\tilde{u}_2 (\sqrt{ ( x_1 - t_1 )^2 + ( x_2 - s_2 )^2} \, , x' ) 
\end{array}
\leqno{(14)}
$$
 for any $ x_1, x_2 \in \R $  and $ x' \in \{ e_1, e_2 \}^{\perp }.$
Since $ u_1 = u_2$ on $ \Pi_{e_1, t_1} = t_1 e_1 + e_1^{\perp }$, it follows that 
$$
\tilde{u}_1 ( | x_2 - t_2 | \, , x' ) 
=
\tilde{u}_2 (| x_2 - s_2  |\, , x' ) 
\leqno{(15)}
$$
for any $ x_2 \in \R$ and  $ x' \in \{ e_1, e_2 \}^{\perp }.$
In particular, (15) implies that for fixed $ x' \in \{ e_1, e_2 \}^{\perp }$,
$ \tilde{u}_1 (  \cdot , x' )$ and $  \tilde{u}_2 (  \cdot , x' )$   are 
periodic of period $ a = 2( s_2 - t_2)$. 
Passing to cylindrical coordinates $ x_1 = t_1 + r \cos \te$, 
$ x_2 = t_2 + r \sin \te$,
 $ x'$ and using Fubini's theorem we have 
$$
\begin{array}{l}
\ds \int_{\Pi_{e_1, t_1 }^-} G_j( u(x), |\nabla u(x) |) \, dx 
= \ds \int_{\Pi_{e_1, t_1 }^-} G_j( u_1(x), |\nabla u_1(x) |) \, dx 
\\
\\
= \ds \int_0^{\infty } \int_{\frac{\pi}{2}}^{\frac{3 \pi}{2}} 
\int_{ \{ e_1, e_2 \}^{\perp} }
G_j ( \tilde{u}_1  ( r, x') , 
|\nabla \tilde{u}_1( r, x') |) \, dx' \, d \te \,r\,  dr
\\
\\
= \pi \ds \int_0^{\infty } \int_{ \{ e_1, e_2 \}^{\perp} }
G_j(  \tilde{u}_1   ( r, x') , 
|\nabla \tilde{u}_1( r, x') |) \, dx' \, r\,  dr.
\end{array}
\leqno{(16)}
$$
Let  $ h_j ( r) = \ds \int_{ \{ e_1, e_2 \}^{\perp} }
G_j ( \tilde{u}_1  ( r, x') , |\nabla \tilde{u}_1( r, x') |) \, dx' $. 
The function  $ h_j$ is well-defined for a.e. $ r \geq 0$,  measurable, 
periodic of period $a$, 
and $ \pi \ds \int_0 ^{\infty } r h_j (r)\, r =  \la _j /2$. 
By periodicity we have 
$ \ds \int_{na }^{(n+1) a} r h_j(r) \, dr = na \int_0^a h_j (r) \, dr + 
 \int_0^a r h_j (r) \, dr$, thus 
$ {\ds  \int_0^{n a}}  r h_j (r) \, dr = \frac{ n( n-1) }{2} a 
\ds  \int_0^a h_j (r) \, dr +  n \int_0^a r h_j (r) \, dr$.
It follows  that necessarily $\ds  \int_0^a h_j (r) \, dr = 0 $ and 
$ \ds \int_0^a r h_j (r) \, dr = 0$ and this implies 
$ \ds \int _0^{\infty } r h_j (r) \, dr = 0$, i.e. $ \la _j = 0$ for any $ j$, 
contrary to the assumptions of Theorem 2. 
Consequently the case 1 b) may never occur. 

\medskip

Case 2:  $ e_2 $ and $ v_2$ are not colinear. 
It is then clear that the space $\mbox{ Span}\{ e_1, e_2, v_2 \}$ is $3-$dimensional 
(thus  $ N \geq 3$). 
Let $ \{e_4, \dots, e_N \}$ be an orthonormal basis of $\{ e_1, e_2, v_2 \}^{\perp }$. 
We choose $ e_3 $ and $ v_3 $ in such a way that 
$\Bo = \{ e_1, \dots, e_N \}$ and $ \Bo '= \{ e_1, v_2, v_3, e_4, \dots, e_N \}$
are orthonormal basis in $ \R^N$ with the same orientation.
There exists $ \te \in ( 0, \pi ) \cup ( \pi, 2 \pi) $ such that 
$ v_2 = \cos \te \, e_2 + \sin \te \, e_3$ and
$ v_3 = - \sin \te \, e_2 + \cos  \te \, e_3$. 
Since $ \sin \te \neq 0$, there exist some $ \al , \beta \in \R$ such that 
$ t_2 e_2 + \al e_3 = s_2 v_2 + \beta v_3$. 
Let $ y = t_1 e_1 + t_2 e_2 + \al e_3$. 
We denote $ u^* = u( \cdot - y)$, 
$ u_1 ^* = u_1 ( \cdot - y)$ and  $ u_2 ^* = u_2 ( \cdot - y)$.
It is obvious that $ u^*, u_1^*$ and $ u_2^*$ are minimizers for ($\Po '$), 
$ u_1 ^*$ is radially symmetric with respect to 
$\mbox{ Span}\{ e_3, \dots, e_N \}$, 
$ u_2 ^*$ is radially symmetric with respect to 
$\mbox{ Span}\{ v_3, e_4, \dots, e_N \}$, 
$ u^* = u_1 ^* $ on $ \Pi_{e_1, 0}^- \cup \Pi_{e_1, 0}$ and 
$ u^* = u_2 ^* $ on $ \Pi_{e_1, 0}^+ \cup \Pi_{e_1, 0}$.
Proceeding as in the proof of Theorem 1 we show that either $u^*$ is radially symmetric with respect to  $\mbox{ Span}\{ e_4, \dots, e_N \}$, 
or there exists $ w \in \mbox{Span}\{ e_2, e_3 \}$, 
such that $u^*$ is radially symmetric with respect to 
$\mbox{ Span}\{ w, e_4, \dots, e_N \}$.
In any case it follows that $u$ is radially symmetric with respect to an affine subspace of dimension at most $ k-1 = N-2$.
This completes the  proof of Theorem 2. 
\hfill $\Box$

\section{Remarks and examples }

{\bf Remark 4. } 
If $\Om $ is connected and a unique continuation principle is available for minimizers, the proofs in the preceding section can be considerably simplified. Moreover, it is possible to deal with $N-1$ constraints in Theorem 1, respectively with $N$ constraints in Theorem 2 (but this is of quite limited interest in applications because we get only symmetry with respect to a hyperplane). 

For example, consider the problem ($\Po 1$) of minimizing  
$$
E(u) = \int_{\Om } \frac 12 | \nabla u |^2 + F(u) \, dx 
\qquad \mbox{ in } 
H^1( \Om, \R^m) \quad (\mbox{or in } H_0^1( \Om, \R^m))
$$
under the constraints $ Q_j(u) = \ds \int _{\Om } G_j(u) \, dx = \la _j$, 
$1 \leq j \leq k$ and the following  standard assumptions:

\medskip

{\bf H1.} $F, G_1, \dots, G_k \in C^2( \R^m, \R)$, $ F(0) = G_j(0) = 0$,  
$ \nabla F (0) =\nabla G_j (0) = 0$, and
$$
|\nabla F(u) | \leq C | u|^p, \qquad
|\nabla G_j (u) | \leq C |u|^p \qquad
\mbox{ for } |u | \geq 1, \mbox{ where } p < \frac{N+2}{N-2}.
$$

\medskip

{\bf H2. } If $ u \in H^1( \Om, \R^m)$ 
(respectively  $ u \in H_0^1( \Om, \R^m)$)
is nonconstant and \\ $ \sum_{j =1}^k \al _ j \nabla G_j(u) = 
\sum_{j =1}^k \beta_ j \nabla G_j(u)$ on $ \Om$, then $ \al _j = \beta _j $ for $ j=1, \dots, k$.

\medskip

Suppose that $u$ is a minimizer for ($\Po 1$) 
and a hyperplane $\Pi $ (with $ 0 \in \Pi $ if $\Om \neq \R^N$) splits the contraints in two for $u$. 
As before, it follows easily that the functions $ u_{\Pi^-}$  and 
$ u_{\Pi^+}$  are  minimizers for ($\Po 1$).
Thus $u$ and $u_{\Pi^-}$   satisfy the Euler-Lagrange equations
$$
- \Delta u + \nabla F(u) + \sum_{j =1}^k \al _ j \nabla G_j(u) = 0
\qquad \mbox{ in } \Om, \quad \mbox{ respectively} 
\leqno{(17)}
$$
$$
- \Delta u_{\Pi ^-} + \nabla F(u_{\Pi ^-}) 
+ \sum_{j =1}^k \beta _ j \nabla G_j(u_{\Pi ^-}) = 0
\qquad \mbox{ in } \Om 
\leqno{(18)}
$$
for some $ \al _1, \dots, \al _k, \beta _1, \dots, \beta _k \in \R$. 
By standard regularity theory we get $ u, \; u_{\Pi^-} \in W^{2,q }(\Om )$ 
for any $ q \in [2, \infty)$. In particular, 
$ u, \; u_{\Pi^-} \in C^{1, \al } (\Om ) $ for $ \al \in [0, 1)$, 
and $  u, \; u_{\Pi^-} $ as well as their derivatives are bounded on $ \Om$. 
If $u$ is constant on $ \Om \cap \Pi ^-$, it follows form (17) and the unique continuation principle (see \cite{lop1}) that $u$ is constant on $ \Om$. 
Otherwise, from (17) and (18) we obtain 
$ \sum_{j =1}^k \al _ j \nabla G_j(u) = 
\sum_{j =1}^k \beta_ j \nabla G_j(u)$ on $ \Om \cap \Pi ^-$ and by {\bf H2 } 
we infer that $\al _j = \beta _j$, $ j =1, \dots, k$. 
Denoting $w$ = $ u - u_{\Pi^-} $,  (17) and (18) imply that $w$ satisfies 
$$
- \Delta w + A(x) w = 0 \qquad \mbox{ in } \Om, 
$$
where $A \in L^{\infty } ( \Om, M_m(\R))$. Since $ w = 0 $ in $ \Om \cap \Pi^-$,
by the unique continuation principle we find $ w = 0$ in $\Om$, 
i.e. $u =  u_{\Pi^-} $ and $u$ is symmetric with respect to $\Pi$. 
Hence we have proved that $u$ is symmetric with respect to any hyperplane that 
splits the constraints in two. 
The rest of the proof is as in the preceding section. 

Note that a nondegeneracy hypothesis like {\bf H2 } is needed to use 
a unique continuation principle.

\bigskip

\noindent
{\bf Remark 5. } 
In Theorems 1 and 2, any supplementary constraint in the minimization problem produces the loss of one direction of symmetry for minimizers. 
Under the general assumptions made there, this loss of symmetry cannot be
avoided, as it can be seen in the following simple examples. 

\medskip

\noindent
{\bf Example 6. } i) Let $ \Om $ be either a ball or an annulus in $ \R^N$, centered at the origin.
Consider $ F, G \in C^2(\R, \R) $ satisfying assumption {\bf H1 } in Remark 4 and such that the problem ($\Po _1 $) 
of  minimizing   
$ E_1 (u) = \int_{\Om } \frac 12 |\nabla u |^2 + F(u) \, dx $ 
in $ H^1(\Om) $ 
 under the constraint $ \int_{\Om }  G(u) \, dx = \la $
admits a nonconstant solution $ u_*$. 
It has been shown in \cite{lop2} that $ u_* $ cannot be radially symmetric 
about $0$ (but, of course, $u_*$ is radially symmetric with respect to a line passing through $0$). 
Consider the problem 
$$
\begin{array}{c} 
\mbox{ minimize } 
E_k(u) = \ds \int_{\Om } \frac 12 |\nabla u |^2 + F(u_1) +\dots + F(u_k)\, dx, 
\\
\\
\mbox{ under the constraints } \ds \int_{\Om }  G(u_j) \, dx = \la, 
\qquad j = 1, \dots, k, 
\end{array}
\leqno{(\Po _k)}
$$
 where $ u = ( u_1, \dots, u_k) \in H^1( \Om , \R^k)$.
It is clear that $ u=( u_1, \dots, u_k)$ is a solution of ($\Po _k $) if and only if each $ u_j$ is a solution of ($\Po _1$). 
If $ R_1, \dots, R_k$ are rotations in $ \R^N$, the function 
$ u(x)  = ( u_* ( R_1 x), \dots, u_*(R_kx))$ is a solution of ($\Po _k $).
We infer that there are minimizers of ($\Po _k $) that are not radially symmetric with respect to any $(k-1)-$dimensional vector subspace of $ \R^N$. 

\medskip

ii) Consider two functions $ F, G \in C^2( \R, \R) $ satisfying assumption 
{\bf H1 } in Remark 4 and $ \la \in \R^*$ such that the problem ($\Po _1 '$) 
consisting in  minimizing 
$\tilde{E}_1 (u) = \int_{\R^N} \frac 12 |\nabla u |^2 + F(u)\, dx $ 
in $ H^1(\R^N) $ 
 under the constraint $ \int_{\R^N}  G(u) \, dx = \la $
admits a nonconstant solution $ \tilde{u} $. 
It folows immediately from 
Theorem 2 that $\tilde{u}$ is radially symmetric with respect to  a point; 
 we may assume that it is radially symmetric about the origin.
It is easy to see that $ u = ( u_1, \dots, u_k) \in H^1( \R^N, \R^k)$ 
is a solution of the problem 
$$
\begin{array}{c} 
\mbox{ minimize } 
\tilde{E}_k(u) = \ds \int_{\R^N } \frac 12 |\nabla u |^2 + F(u_1) +\dots + F(u_k)\, dx, 
\\
\\
\mbox{ under the constraints } \ds \int_{\R^N }  G(u_j) \, dx = \la, 
\qquad j = 1, \dots, k, 
\end{array}
\leqno{({\Po } _k ' )}
$$
in $H^1( \R^N, \R^k)$  
if and only if each $ u_j $  is a solution of ($\Po _1 '$). 
Therefore for any $ y_1, \dots , y_k \in \R^N$, the function 
$ u = ( u_1 ( \cdot + y_1), \dots, u_k( \cdot + y_k))$ 
is a solution for ($ {\Po } _k '$). 
Obviously, this minimizer is radially symmetric with respect to some $(k-1)-$dimensional affine subspace but, in general, it is not radially symmetric with respect to any affine subspace of lower dimension.

In Example 6, the loss of symmetry comes from the fact that problems 
($ {\Po } _k $)
and ($ {\Po } _k '$) 
are decoupled: they can be decomposed into $k$ independent scalar problems, each
 of them being rotation (respectively translation) invariant.
It is then natural to ask whether in 
general problems like ($\Po $) or ($\Po '$) 
the loss of directions of symmetry could exceed the number 
of components of minimizers.
The answer is affirmative, as it can be seen in the next example 
which shows that,   
in general, the result of Theorem 2 is optimal 
even for scalar-valued minimizers. 

\medskip
\noindent
{\bf Example 7. } 
We construct here a  minimization problem of the form ($\Po '$) 
involving two constraints and whose real-valued minimizers are {\it  not} radial
with respect to a point (of course, these minimizers are axially symmetric). 
This example relies on the existence of a nonnegative minimizer with compact support for 
a problem involving one constraint. A similar construction has already been used in 
\cite{brock-JMAA}.

Let $ f \in C(\R) \cap C^1(0, \infty )$ be a real-valued function
satisfying the following conditions:

\medskip

\noindent
{\bf C1. } $f(s) = 0 $ on $(-\infty, 0]$ and $ f(s) = s^{\al }$ for $ s \in (0, 1]$, 
where $ \al \in (0, 1)$.

\medskip

\noindent
{\bf C2. } The function $ F(s) := \int_0^s f(\tau )\, d \tau $ has compact support.

\medskip

\noindent
{\bf C3. } There exists $ \zeta > 0$ such that $ F(\zeta ) < 0$. 

\medskip

Let $ N \geq 3$ and $ \Xo = \DR \cap L^{1 + \al } ( \R^N)$. 
We introduce the functionals $ T(u) = \ds \int_{\R^N} |\nabla u |^2   \, dx $ and 
$V(u) = \ds \int_{\R^N} F(u(x)) \, dx$. 
It is clear that $ F(u) \in L^1(\R^N)$ for any $ u \in \Xo$ and $T, \; V$ are 
well-defined, $C^1$ functionals on $ \Xo$. 
We consider the minimization problem:
$$
\mbox{ minimize } T(u) \mbox{ in } \Xo \mbox{ subject to the constraint } V(u) = -1.
\leqno{(\Mo _1)}
$$

We denote $ I = \inf \{T(u) \; | \; u \in \Xo, V(u) = -1 \}$ 
and we proceed in several steps. 

\medskip

{\it Step 1.} {\it We have $ I >0$ and problem $(\Mo _1) $ has a minimizer  $u_* \in \Xo $.}
The proof of this fact is a straightforward modification of the proof of 
Theorem 2  in  \cite{BL}
or of the proof of Theorem 1 in \cite{ferrero-gazzola}, so we omit it. 

\medskip

{\it Step 2.} {\it Any minimizer $u$ of $ (\Mo _1)$ is nonnegative, bounded, $C^1$, 
has compact support and satisfies the 
equation   
$
- \Delta u + \beta _0 f(u) = 0 $ in $ {\mathcal{D}}'(\R^N)$, 
 where $ \beta _0 = \frac{N-2}{2N}I$. }

Let $ u^+ = \max (u, 0)$ and $ u^- = \max (-u, 0)$. Then $ u^+, u^- \in \Xo$, 
$V(u^+) = V(u) = -1 $ and $T(u) = T(u^+) + T(u^-) \geq T(u^+)$. 
Since $u$ is a minimizer, we must have $T(u^+) = T(u) $ and $ T(u^-) = 0$, hence 
$ u^- = 0 $ in $ \DR$, that is $ u \geq 0 $ a.e. 
Take $ C>0$ such that $\mbox{supp}(F) \subset [0, C]$ and denote $ u_0 = \min (u, C)$, 
$u_C = \max ( u-C, 0)$. It is obvious that 
$ u_0, u_C \in \Xo $, $ u = u_0 + u_C$, $V(u_0) = V(u) = -1 $ and 
$ T(u) = T(u_0) + T(u_C)$.
As above we infer that $ T(u_C) = 0$, consequently $ u_C = 0 $ in $\DR$ and 
$ u \leq C $ a.e.

Since $T$ and $V$ are $C^1$ functionals on $\Xo$, it is easy to see that 
$u$ satisfies an Euler-Lagrange equation $T'(u) + 2 \beta V'(u) = 0 $ in $\Xo '$
for some $ \beta \in \R$ and this implies
$$
- \Delta u + \beta f(u) = 0 \qquad \mbox{ in } {\mathcal{D}}'(\R^N).
\leqno{(19)}
$$
Since $u \in L^{\infty } (\R^N)$ and $f$ is continuous, 
by standard elliptic estimates it follows that 
$ u \in W_{loc}^{2, p} (\R^N) $ for any $ p \in (1, \infty)$, 
thus $ u \in C_{loc}^{1, \g }(\R^N) $ for $ \g \in [0, 1)$. 
In particular,  $u$ is $C^1$.

It is standard to prove that $u$ satisfies the Pohozaev identity 
$(N-2) T(u) + 2 \beta N V(u) = 0$ 
(to see this, it suffices to multiply (19) by $\chi ( \frac xn ) \sum_{i =1 }^N 
x_i \frac{ \p u}{\p x_i} $, where $\chi \in C_c^{\infty }(\R^N) $ 
is such that $ \chi \equiv 1$ on $B(0, 1)$, 
to integrate by parts and then to pass to the limit as $ n \lra \infty $). 
Since $ V(u) = -1 $ and $T(u) = I$, we find $ \beta = \frac{N-2}{2N} I = \beta _0 >0$. 

Let $ v(x) = u(\frac{x}{\sqrt{ \beta _0}})$. Then $ v \in C^1(\R^N)$, $ v \geq 0 $ and $v$ 
satisfies the equation 
$$
 - \Delta v + f(v) = 0 \qquad \mbox{ in } {\mathcal{D}}'(\R^N).
$$
Moreover, we have 
$\ds \int_0^1 \frac{1}{(F(s))^{\frac 12}} \, ds 
= (\al + 1) ^{\frac 12} \int_0^1 \frac{1}{s^{\frac{\al +1}{2}}} \, ds 
< \infty$. 
Thus we may use Theorem 2 p. 773 in \cite{PSZ} and we infer that $v$ 
has compact support.
Hence $u$ has compact support. 

\medskip

{\it Step 3.} {\it Any minimizer $u$ of $ (\Mo _1)$ is radially symmetric with respect to 
a point.}
Indeed,  steps 1 and 2 show that $(\Mo _1)$ satisfies assumptions {\bf A1'}
and {\bf A2} in Introduction, 
hence the radial symmetry of minimizers follows from Theorem 2. 
Note that the unique continuation principle is not valid for minimizers 
of $(\Mo _1)$, therefore the method in \cite{lop1} 
cannot be used to prove their radial symmetry. 

\medskip

{\it Step 4.} {\it Construction of nonradial minimizers for a minimization problem 
involving two constraints.}

We introduce the functional $W(u) = \ds \int_{\R^N} F(-u(x)) \, dx$. 
Clearly, $W$ is well-defined and $C^1$ on $ \Xo $. 
We consider the minimization problem:
$$
\mbox{ minimize } T(u) \mbox{ in } \Xo \mbox{ subject to the constraints } V(u) = -1
\mbox{ and } W(u) = -1.
\leqno{(\Mo _2)}
$$

We claim that $ u \in \Xo $ is a solution of $ (\Mo _2)$ if and only if 
$u^+$ and $ u^-$ are solutions of $(\Mo _1)$. 

To see this, let $u_*$ be a minimizer of $(\Mo _1)$, 
radially symmetric with respect to the origin. 
Let $R > 0$ be such that $\mbox{supp} (u_*) \in B(0, R)$. 
For $ y \in \R^N \setminus B(0, 2R)$, we put $u_y (x) = u_* (x) - u_*(x +y) $. 
It is obvious that $ V(u_y) = V(u_*) =-1$, 
$W(u_y ) = V(u_* ( \cdot + y)) = -1 $ and $T(u_y) = T( u_*) + T( u_*(\cdot + y)) = 2I$. 

For any $ u \in \Xo $ satisfying $ V(u) = W(u) = -1$ we have 
$V(u^+) = V(u ) = -1$  and $ V(u^-) = W(u) = -1$,
hence $ T(u^+) \geq I$ and $ T(u^-) \geq I$, consequently $T(u) \geq 2I$. 
We conclude that for any $|y| \geq 2R$, $u_y$ is a minimizer of $(\Mo _2)$.
Moreover, a function $ u \in \Xo $ can solve $(\Mo _2)$ if and only if 
$V(u^+) = V(u^-) = -1$ and $T(u^+) = T(u^-) = I$, i.e. if and only if 
$ u^+$ and $ u^-$ solve $(\Mo _1)$.

As in step 2 we infer that all minimizers of $(\Mo _2)$ are $C^1$. 
Thus $(\Mo _2)$ satisfies the assumptions {\bf A1'} and {\bf A2} 
and Theorem 2 implies that all minimizers of $(\Mo _2)$ are axially symmetric.
Since $ u_* $ is radial with respect to the origin, 
it is clear that any of the minimizers 
$u_y$ is axially symmetric with respect to the line $Oy$, 
but is not radial about a point. 
Hence $(\Mo _2)$ admits nonradial minimizers.

In fact, with some extra work it can be proved that the suport of any minimizer of 
$(\Mo _1)$ is precisely a ball. If  $u$ is a minimizer of $(\Mo _2)$, 
$\mbox{supp} (u) = \mbox{supp}(u^+) \cup \mbox{supp}(u^-)$ is the union of two balls 
with disjoint interiors. 
Therefore no minimizer of  $(\Mo _2)$ can be radially symmetric.

\medskip

In some particular cases, however, minimizers may have more symmetry than 
provided by Theorems 1 and 2, as it can be seen in the following example.

\medskip

\noindent
{\bf Example 8. } Consider the problem of minimizing 
$ E(u) = {\ds \int_{\R }}  \frac 12 |u'(x)| ^2 + F(u(x))\,  dx $ 
in $H^1(\R)$, 
under an arbitrary number of constraints 
$ {\ds \int_{\R}}G_j(u(x)) \, dx = \la _j$, $ 1 \leq j \leq k$. 
We assume that the functions  $F, G_1, \dots , G_j$ satisfy 
the assumption {\bf H1 }in Remark 4.

In this case Theorem 2 gives no information about the minimizers. 
However, if the problem above admits minimizers, any of them must be symmetric with respect to a point. 
Indeed, let $u$ be a nonconstant minimizer. 
Then it satisfies an Euler-Lagrange equation
$$
- u'' + F'(u) + \al _1 G_1 '(u) + \dots + \al _k G_k'(u) = 0 
\qquad \mbox{ in } \R. 
\leqno{(20)}
$$
It follows easily from (20) that $ u \in C^2(\R, \R)$. 
Since $ u(x) \lra 0  $ as $ x  \lra  \pm \infty$, 
$u$ achieves its maximum or its minimum at some point $a \in \R$ and consequently $ u'(a) = 0$. 
Let $ \tilde{u}(x) = u( 2a -x)$.
Then $\tilde{u}$ satisfies (20) and $ \tilde{u}(a) = u(a)$, 
$ \tilde{u}'(a) = u'(a)=0$.
Since the Cauchy problem associated to (20) has unique solution, we have 
$ u = \tilde{u}$, i.e. $u$ is symmetric about $a$. 
Moreover, we see that $u$ must be symmetric with respect to any of its critical points. 
Since $u$ cannot be periodic, we infer that there are no other critical points, 
thus $u$ is monotonic on $(-\infty, a]$ and on $[a, \infty )$. 

\medskip

We have discussed in the first section an example of problem where arbitrarily many  constraints were allowed and the symmetry properties of minimizers 
did not depend on the number of constraints (see \cite{brock-JMAA}). This fact is due to the assumptions
made on the nonlinear term (monotonicity in $|x|$ and cooperativity condition),
that imply a strong coupling between the components of the minimizers and 
prevent situations like those in Examples 6 and 7 to occur.

\medskip

\noindent
{\bf Remark 9. } Our results can be extended in an obvious way to minimization problems on cylinders.
To be more specific, consider the problem ($\Po _c $)  consisting in  minimizing
$$
E(u) = \int _A \int_{\Om  } 
F(|x|, y, u(x,y), |\nabla _x u (x, y)|, \nabla _y u(x,y), \dots, 
\nabla _y ^{\ell } (x, y) ) \, dxdy 
$$
under the constraints 
$$
Q_j(u) = \int_A \int_{\Om  } 
G_j(|x|, y, u(x,y), |\nabla _x u (x, y)|, \nabla _y u(x,y), \dots, 
\nabla _y ^{\ell } (x, y) ) \, dxdy , 
\quad j =1, \dots, k,
$$
where $ x \in \Om \subset \R^{N_1}$, $ y \in A \subset \R^{N_2}$, 
$ \Om $ is an open set invariant by rotations in $ \R^{N_1}$ and $A$ is a measurable set in $ \R^{N_2}$. 
We assume that problem ($\Po _c $)
admits 
 minimizers in a functional space $ \Xo $ and  the following assumptions hold: 

\medskip

{\bf {A1}$_c$.} For any $ w \in \Xo $ and any hyperplane $ \Pi $ in $ \R^{N_1}$ containing the origin, we have 
$w_{(\Pi \times \R^{N_2}) ^-}, w_{(\Pi \times \R^{N_2}) ^+} \in \Xo$. 

\medskip

{\bf {A2}$_c$.} For any minimizer $ u \in \Xo $ and any $ y \in A$, 
the function $u ( \cdot, y)$ is $C^1$ on $\Om$. 

\medskip

Note that the minimization problem may involve derivatives 
of any order
in $y$ and 
we do not need more regularity of minimizers with respect to $y$
than provided by the fact that $ u \in \Xo$.  

We have the following results, the proofs being similar to those of Theorems 1 and 2.

\bigskip

\noindent
{\bf Theorem 1'. } \it Assume that $ u$ is a minimizer for problem ($\Po _c $)
in $ \Xo$, assumptions {\bf {A1}$_c$} and {\bf {A2}$_c$} are satisfied and 
$ 0 \leq k \leq N_1-2$. 
There exists a $k-$dimensional vector subspace $V$ of $ \R^{N_1}$ such that $u$ is radially symmetric with respect to 
$ V \times \R^{N_2}$. 

\rm 

\bigskip

\noindent
{\bf Theorem 2'. } \it Assume that $ \Om = \R^{N_1}$, $ 1 \leq k \leq N_1-1$ and the functions $F$, $G_j $ in ($\Po _c $) do not depend on $x$. 
Assume also that {\bf {A2}$_c$} is satisfied  and {\bf {A1}$_c$} 
holds for any affine hyperplane $ \Pi $ in $ \R^{N_1}$.
If $ u $ is a minimizer for problem ($\Po _c $) in $\Xo$, 
there exists a $(k-1)-$dimensional affine subspace $ V \subset \R^{N_1}$ 
such that 
$u$ is radially symmetric with respect to $ V \times \R^{N_2}$. 

\rm

\bigskip

{\it Acknowledgements. } I am very grateful to Petru Mironescu, Petru Jebelean and Alberto Farina for interesting and helpful discussions.

\end{document}